\documentclass[12pt]{amsart}
\usepackage{amstext,amsmath,amssymb,amscd,amsthm,amsfonts,amstext,amsbsy,graphicx,latexsym,multirow}
\usepackage{latexsym}
\usepackage[english]{babel}
\usepackage[latin1]{inputenc}
\newtheorem{defis}{Definition}[section]
\newtheorem{theos}[defis]{Theorem}
\newtheorem{lemm}[defis]{Lemma}
\newtheorem{prop}[defis]{Proposition}
\newtheorem{corol}[defis]{Corollary}

\newtheorem{problem}[defis]{Problem}
\newtheorem{example}[defis]{Example}

\newtheorem{fact}{Fact}
\newenvironment{demos}{\medskip\noindent\textsc{Proof.}}{\medskip}

\title{Are Eberlein-Grothendieck scattered spaces $\sigma$-discrete?}
\author{Antonio Avil\'{e}s and David Guerrero S\'{a}nchez}

\address{Departamento de Matem\'{a}ticas, Universidad de Murcia, 30100 Murcia, Spain}
\email{avileslo@um.es,david.guerrero@um.es}

\begin{document}
\maketitle

\begin{center}
\textbf{Abstract}
\end{center}
{\small A space $X$ is Eberlein-Grothendieck if $X\subset C_p(K)$
for some compact space $K.$ In this paper we address the problem of
whether such a space $X$ is $\sigma$-discrete whenever it is
scattered. We show that if $w(K)\leq\omega_1$ then $X$ is
$\sigma$-dicrete whenever $X$ has height $\omega_1$ and it is
locally compact or locally countable. It is also proved that every
Lindelöf \v{C}ech-complete scattered space is $\sigma$-compact.}\\


\noindent {\small \textbf{Mathematical subject code}: 54C35.}

\section{Introduction}

In [Ar] Arhangel'skii defined Eberlein-Grothendieck spaces as those
homeomorphic to a subspace of $C_p(K)$ for some compact space $K.$
Notice that if $X$ is a subset of a Banach space $E$ with the weak
topology then $X$ embeds in $C_p(\textbf{B}E^*)$ hence $X$ is
Eberlein-Grothendieck. The main purpose of this paper is to
introduce the following problem:
\begin{problem}\label{mainproblem}
Are
Eberlein-Grothendieck scattered spaces $\sigma$-discrete?
\end{problem}
 A
particular case of this problem was posed in [Hay] where Haydon
asked if for every compact $K$ the space $C_p(K,\{0,1\})$ is
$\sigma$-discrete whenever it is scattered. This kind of question is
related to the following notions introduced by J.E. Jayne, I.
Namioka and C.A. Rogers in [JNR]. Given a set $X$, a metric $\rho$
on $X$ and $\varepsilon>0,$ a family $\mathcal{A}$ of subsets of $X$
is $\varepsilon$-small if $diam_\rho(A)<\varepsilon$ for every
$A\in\mathcal{A}.$ A topological space $(X,\tau)$ has the property
SLD with respect to a metric $\rho$ on $X$ if for every
$\varepsilon>0$ there is a countable cover $\{X_n:n\in\omega\}$ of
$X$ such that for each $n\in\omega$ the space $X_n$ admits a
$\tau$-open cover which is $\varepsilon$-small. On the other hand a
topological space $(X,\tau)$ is $\sigma$-fragmented by a metric
$\rho$ on $X$ if for every $\varepsilon>0$ there is a countable
cover $\{X_n:n\in\omega\}$ of $X$ such that for each $n\in\omega$
and every $Y\subset X_n$ there exists a nonempty relative
$\tau$-open subset $U$ of $Y$ with $diam_\rho(U)<\varepsilon.$ It is
clear that if a topological space has SLD with respect to some
metric, then it is $\sigma$-fragmented as well, but the following is
an open question:

\begin{problem}\label{SLDproblem}
Are the properties of $\sigma$-fragmentability and SLD equivalent
when $X$ is a Banach space endowed with its weak topology and with
its norm metric, or when $X$ is of the form $C_p(K)$ endowed with
the uniform metric?
\end{problem}

This problem has its origin in the theory of renorming of Banach
spaces, as these properties are conjectured as possible internal
characterizations of Banach spaces admitting a norm with the
Kadets-Klee property. We refer to [MOTV, Section 3.2, p.54] for
information about this topic. It is easy to see that if a space with
the discrete metric is SLD then it is $\sigma$-discrete and if it is
$\sigma$-fragmented then it is $\sigma$-scattered. Thus Problem
\ref{mainproblem} can be viewed as the discrete version of Problem
\ref{SLDproblem}. Indeed the following question is open as well:

\begin{problem}\label{CK2problem}
If $C_p(K,\{0,1\})$ is scattered (respectively $\sigma$-discrete),
does it imply that $C_p(K)$ is $\sigma$-fragmentable (respectively
SLD)?
\end{problem}

A version of this problem states the same thing considering the weak
topology of $C(K)$ instead of the pointwise convergence topology.
The answer is known to be positive when $K$ is scattered (see [Hay]
and [Mtz].) Since the restriction of the uniform metric of $C_p(K)$
to $C_p(K,\{0,1\})$ is discrete, a positive answer to Problems
\ref{mainproblem} and \ref{CK2problem} combined would give a
positive answer to
Problem \ref{SLDproblem} in the case of spaces of continuous functions.\\

We observe that it follows from known results that Problem
\ref{mainproblem} has positive solution when $X$ is compact: Alster
proved in [Al] that an Eberlein compact scattered space $K$ is
strong Eberlein which implies that $K$ embeds into $\{0,1\}^\Gamma$
for some $\Gamma$ and $|supp(x)|<\omega$ for every $x\in K.$ For
each $n\in\omega$ we can define $X_n=\{x\in K: |supp(x)|=n\}$ hence
we can write $K=\bigcup X_n$ where each $X_n$ is discrete. In this
paper we will prove some generalizations of this fact.

In Section 3, we give some partial positive answers to Problem
\ref{mainproblem}. The first one is the following:

\begin{theos}\label{locallycompacthm}
If $X$ is an Eberlein-Grothendieck locally compact scattered space
of height lower than $\omega_1\cdot\omega_1$, then $X$ is
$\sigma$-discrete.
\end{theos}

In the same section we will prove also our second result which
states:

\begin{theos}\label{locallycountablethm}
If $X$ is an Eberlein-Grothendieck locally countable scattered space
of cardinality $\omega_1$, then $X$ is $\sigma$-discrete.
\end{theos}

Remember that a transfinite sequence $\{x_\alpha : \alpha<\lambda\}$
of elements of a topological space is right-separated if for every
$\mu<\lambda$ there is an open set $U$ such that $U\cap \{x_\alpha :
\alpha<\lambda\} = \{x_\alpha : \alpha<\mu\}$. A topological space
is scattered if and only if it can be written as a right separated
sequence $X = \{x_\alpha : \alpha<\lambda\}$. From this point of
view, Theorem \ref{locallycountablethm} implies that Problem
\ref{mainproblem} has a positive solution in the first non-trivial
case, when $\lambda=\omega_1.$

\begin{corol}
If $X = \{x_\alpha : \alpha<\omega_1\}\subset C_p(K)$ is a
right-separated $\omega_1$-sequence, then $X$ is $\sigma$-discrete.
\end{corol}

In both of our results mentioned above, $X$ is homeomorphic to some
$X'\subset C_p(K)$ where $K$ has weight $\omega_1$. By [DJP, Theorem
1.2] in this case the space $X$ is hereditarily meta-Lindelöf. This
is the hypothesis that we shall actually assume, so that Theorem
\ref{locallycompacthm} is proved by applying the ideas developed in
[Al] to show that every open cover of a hereditarily meta-Lindelöf
locally compact scattered space of height at most $\omega_1$ has a
point finite clopen refinement, while Theorem
\ref{locallycountablethm} is proved by showing that hereditarily
meta-Lindelöf locally countable scattered spaces are
$\sigma$-discrete. The latter fact is stated in [HP] without proof,
but the argument that they suggest does not seem to be correct.
Section 3 ends with a corollary that, at least when $K$ is
scattered, the SLD property of $C_p(K)$ can be characterized as a
kind of
$\omega_1$-$\sigma$-fragmentability.\\

In Section 4 we deal with another generalization of the Eberlein
compact scattered spaces which is the class of Eberlein-Grothendieck
Lindelöf $\Sigma$ scattered spaces. In this case we show that these
spaces are $\sigma$-discrete as an easy consequence of the results
in [Ha] and [Ny]. But for the special subcase of the
Eberlein-Grothendieck Lindelöf \v{C}ech-complete scattered spaces we
actually prove their $\sigma$-compactness applying the methods of
topological games developed in [Te].\\

Since the results in Section 3 depend very strongly on the property
of hereditary metalindelöfness of the Eberlein-Grothendieck spaces
considered in that section, we decided to ask if such property is
enough for a scattered topological space to be $\sigma$-discrete. In
other words, is every hereditarily meta-Lindelöf scattered space
$\sigma$-discrete? We show in Section 5 that this is not the case in
general by constructing a hereditarily meta-lindelöf scattered space
that is not $\sigma$-discrete. However, the space constructed in
section 5 is not even a Hausdorff space, therefore it is not
Eberlein-Grothendieck.\\

In section 6 we pose some of the remaining open problems related to
determine whether Eberlein-Grothendieck scattered spaces are
$\sigma$-discrete. The general question here is still open so we
consider several subcases. The most important of these subcases is
stated in Problem 6.1 which reads: is every Eberlein-Grothendieck
Lindelöf scattered space $\sigma$-discrete?

\section{Notation and terminology}

Unless otherwise stated, every topological space in this article is
assumed to be Tychonoff. The topology of $X$ is denoted by $\tau(X)$
and $\tau^*(X)$ is the family of non-empty open subsets of $X.$ For
$C\subset X$ the family of all open sets of $X$ that contain $C$ is
denoted by $\tau(C,X);$ if $x\in X$ then we write $\tau(x,X)$
instead of $\tau(\{x\},X).$

The space of all continuous functions from a space $X$ into a space
$Y$, endowed with the topology inherited from the product space
$Y^X,$ is denoted by $C_p(X,Y).$ The space $C_p(X,\mathbb{R})$ will
be abbreviated by $C_p(X).$ For every $f\in C_p(X,Y),$ define the
dual map $f^*:C_p(Y)\rightarrow C_p(X)$ by $f^*(g)=g\circ f$ for
every $g\in C_p(Y).$ For a space $X\subset C_p(Y)$ the diagonal
product of the elements of $X$ is the map $\varphi:Y\rightarrow
\mathbb{R}^X$ defined by $\varphi(y)(f)=f(y)$ for any $y\in Y$ and
$f\in X.$

A space $X$ is Eberlein-Grothendieck if there is a compact space $K$
such that $X$ is homeomorphic to a subspace of $C_p(K).$ A space is
scattered if every subspace of it contains an isolated point. Given
a family $\mathcal{F}$ of subsets of a space $X,$ a refinement of
$\mathcal{F}$ is a family $\mathcal{G}$ such that
$\bigcup\mathcal{F}=\bigcup\mathcal{G}$ and for each
$G\in\mathcal{G}$ there is $F\in\mathcal{F}$ such that $G\subset F.$
A $\sigma$-compact ($\sigma$-countably compact) space is the
countable union of compact (countably compact) spaces. A space is
cosmic if it has a countable network.

As in [Bu], a space $X$ is metacompact if every open cover of $X$
has a point finite open refinement. Given a family $\mathcal{F}$ of
subsets of a space $X,$ for every point $x\in X,$ define
$ord(x,\mathcal{F})=|\{F\in\mathcal{F}: x\in F\}|.$ A space $X$ is
weakly $\theta$-refinable if every open cover of $X$ has an open
refinement
$\mathcal{V}=\displaystyle\bigcup_{n\in\omega}\mathcal{V}_n$ such
that for each $x\in X$ there is $n\in\mathbb{N}$ such that $1\leq
ord(x,\mathcal{V}_n)<\omega.$ The rest of our terminology is
standard and can be found in the books [Ar], [Tk2] and [En].

\section{Hereditarily meta-Lindelöf Eberlein-Grothendieck spaces}

The first mandatory step is to pose the problem as one of covering
properties; we can do it by recalling that Nyikos established in
[Ny, Theorem 3.4] that a scattered space is $\sigma$-discrete if and
only if it is hereditarily weakly $\theta$-refinable. It is clear
that this result implies that every metrizable scattered space is
$\sigma$-discrete as it was proved by Telgarsky in [Te, Theorem
12.1] by other methods. Yakovlev showed in [Ya] that every Eberlein
compact space is hereditarily $\sigma$-metacompact hence
hereditarily weakly $\theta$-refinable, so it suffices to invoke
[Ar, Corollary III.4.2] to conclude that every Eberlein-Grothendieck
pseudocompact scattered space is $\sigma$-discrete.

In [DJP, Teorem 1.2] it was proved by Dow, Junnila and Pelant that
if $K$ is a compact space of weight $\omega_1$ then $C_p(K)$ is
hereditarily meta-Lindelöf. It follows that every
Eberlein-Grothendieck space of cardinality $\omega_1$ is
hereditarily meta-Lindelöf. Indeed, suppose that $X\subset C_p(K)$
where $K$ is a compact space such that $w(k)>\omega_1=|X|.$ Let
$\varphi(p)(f)=f(p)$ for any $p\in K$ and $f\in X.$ In [Tk2, Problem
166] it is proved that $\varphi(p)\in
C_p(X)\subset\mathbb{R}^{\omega_1},$ for any $p\in K$ and that the
map $\varphi:K\rightarrow C_p(X)$ is continuous. The map $\varphi$
thus defined is called the diagonal product of the elements of $X.$
Let $L=\varphi(K)\subset\mathbb{R}^{\omega_1},$ it is clear that
$w(L)=\omega_1.$ Besides, the function $\varphi^*:C_p(L)\rightarrow
C_p(K)$ defined by $\varphi^*(f)= f\circ\varphi$ embeds $C_p(L)$
into $C_p(K)$ as a closed subspace (see [Tk2, Problem 163]). For
every $g\in X,$ if $\pi_g:L\rightarrow \mathbb{R}$ is the projection
of $L$ onto the factor determined by $g,$ then $g=\varphi^*(\pi_g).$
It follows that the space $X$ embeds into $C_p(L)$ which is
hereditarily meta-Lindelöf. We can conclude that every
Eberlein-Grothendiek space of cardinality $\omega_1$ is hereditarily
meta-Lindelöf.

In this section we will strongly exploit metalindelöfness in
Eberlein-Grothendieck spaces combined with a fundamental fact proved
by Alster in [Al] which establishes that every locally countable
family of compact scattered open subspaces of a space $X$ has a
point finite clopen refinement. In order to apply these results to
the study of locally compact scattered spaces we will apply some of
the ideas introduced in [Sp].

\begin{lemm}
Every open cover of a hereditarily meta-Lindelöf locally compact
scattered space of countable height has a point finite clopen
refinement.
\end{lemm}

\begin{demos}
Let $\kappa<\omega_1$ be the height of $X.$ For each $\alpha<\kappa$
denote by $X^{(\alpha)}$ the $\alpha$-th scattering level of the
space $X.$ Recalling that a scattered locally compact space is zero
dimensional, it is easy to see that every open cover $\mathcal{O}$
of $X$ has a clopen refinement $\mathcal{U}$ such that for every
$\alpha <\kappa$ and each $x\in X^{(\alpha)}$ there is
$U_x\in\mathcal{U}$ such that $U_x\cap X^{(\alpha)}=\{x\}$ and
$U_x\cap X^{(\beta)}=\varnothing$ for every $\beta>\alpha.$ We will
show that $\mathcal{U}$ has a point finite clopen refinement.

Fix $\alpha<\kappa.$ The space $Y_\alpha=\bigcup\{U_x:x\in
X^{(\alpha)}\}$ is meta-Lindelöf therefore the cover $\{U_x:x\in
X^{(\alpha)}\}$ of $Y_\alpha$ has a point countable refinement
$\mathcal{U}_{\alpha}.$ For each $x\in X^{(\alpha)}$ we can choose a
single $U^\alpha_x\in\mathcal{U}_\alpha$ such that $x\in
U^\alpha_x\subset U_x.$ Moreover, for each $x\in X^{(\alpha)}$ there
is a clopen set $V^\alpha_x$ such that $x\in V^\alpha_x\subset
U^\alpha_x.$ Observe that if $x,y\in X^{(\alpha)}$ and $x\neq y$
then $V^\alpha_x\neq V^\alpha_y$ which implies $U^\alpha_x\neq
U^\alpha_y$ because $U^\alpha_x\cap
U^{(\alpha)}=\{x\}=V^\alpha_x\cap U^{(\alpha)}.$ Therefore, since
the family $\mathcal{U}_\alpha$ is point countable so is the family
$\{V^\alpha_x:x\in X^{(\alpha)}\}.$ We can now apply [Al,
Proposition] to find a point finite clopen refinement
$\mathcal{V}_\alpha$ of the family $\{V^\alpha_x:x\in
X^{(\alpha)}\}.$

Finally, the family
$\mathcal{V}=\displaystyle\bigcup_{\alpha\in\kappa}\mathcal{V}_\alpha$
is a point countable clopen refinement of $\mathcal{U}.$ Apply [Al,
Proposition] once more to conclude that $\mathcal{V},$ and hence
$\mathcal{U}$, has a point finite clopen refinement.

\end{demos}

\begin{theos}
Every open cover of a locally compact scattered hereditarily
meta-Lindelöf space of height $\omega_1$ has a point finite clopen
refinement.
\end{theos}

\begin{demos}
Let $\mathcal{O}$ be an open cover of $X.$ For each
$\alpha<\omega_1$ denote by $X^{(\alpha)}$ the $\alpha$-th
scattering level of $X.$ As in Lemma 3.1 the open cover
$\mathcal{O}$ of $X$ has a clopen refinement $\mathcal{U}'$ such
that for every $\alpha <\omega_1$ and each $x\in X^{(\alpha)}$ there
is $U_x\in\mathcal{U}'$ such that $U_x\cap X^{(\alpha)}=\{x\}$ and
$U_x\cap X^{(\beta)}=\varnothing$ for every $\beta>\alpha.$ We will
show that $\mathcal{U}'$ has a point finite clopen refinement.

There is a point countable open refinement $\mathcal{U}$ of the
cover $\mathcal{U}'.$ We will show that $\mathcal{U}$ has a point
finite clopen refinement. Fix $U\in\mathcal{U},$ let
$\kappa<\omega_1$ be the height of $U.$ For each $\alpha<\kappa$ and
every $x\in U^{(\alpha)}$ there is a clopen $V^\alpha_x\subset U$
such that $V_x\cap U^{(\alpha)}=\{x\}$ and $V^\alpha_x\cap
U^{(\beta)}=\varnothing$ for every $\beta>\alpha.$ Let
$\mathcal{V}_\alpha=\{V^\alpha_x: x\in U^{(\alpha)}\}.$ The space
$U$ is a hereditarily meta-Lindelöf locally compact scattered space
of countable height, therefore we can invoke Lemma 3.1 to find a
point finite clopen refinement $\mathcal{C}_U$ of the family
$\mathcal{V}=\displaystyle\bigcup_{\alpha\in\kappa}\mathcal{V}_\alpha.$
It follows that the family
$\mathcal{C}=\displaystyle\bigcup_{U\in\mathcal{U}}\mathcal{C}_U$ is
a point countable clopen refinement of $\mathcal{U}$ and its
elements are compact scattered subspaces of $X.$ By [Al,
Proposition] the cover $\mathcal{C},$ and hence $\mathcal{U}$ has a
point finite clopen refinement.

\end{demos}

\begin{corol}
Every hereditarily meta-Lindelöf locally compact scattered space of
height $\omega_1$ is metacompact.
\end{corol}

In the proof of [Ny, Theorem 3.4] Nyikos proved implicitly the
following fact.

\begin{fact}

Let $\kappa$ be a limit ordinal and suppose that $X$ is a scattered
space of height $\kappa$ for which every $U\in\tau(X)$ of height
less than $\kappa$ is $\sigma$-discrete. The space $X$ is
$\sigma$-discrete if it is weakly $\theta$-refinable.

\end{fact}

Metacompactness implies weakly $\theta$-refinability, thus we obtain
the following consequence.

\begin{theos}

Every hereditarily meta-Lindelöf locally compact space $X$ of height
lower than $\omega_1\cdot\omega_1$ is $\sigma$-discrete.

\end{theos}

\begin{demos}
Take any $U\in\tau(X),$ let $\kappa$ be the height of $U$ and
suppose that $\omega_1\leq\kappa<\omega_1\cdot\omega_1.$ Then
$\kappa=\omega_1\cdot\alpha+\beta$ for some countable ordinals
$\alpha,\beta.$ It is easy to see that Theorem 3.2 implies that the
set $\bigcup\{U^{(\gamma)}: \gamma<\omega_1\cdot\alpha\}$ is
$\sigma$-discrete. It follows that $U$ is $\sigma$-discrete and we
can apply the metacompactness of $X$ and Fact 1 to conclude that $X$
is $\sigma$-discrete.
\end{demos}

\begin{theos}

Every scattered hereditarily meta-Lindelöf locally countable space
is $\sigma$-discrete.

\end{theos}

\begin{demos}
Let $\mathcal{C}$ be a point countable cover of the space $X$ such
that every $C\in\mathcal{C}$ is countable. For each $x\in X,$ define
$U^x_0=\bigcup\{C\in\mathcal{C}: x\in C\}.$ Suppose that for every
$x\in X$ we have defined the countable open set $U^x_n\in\tau(x,X).$
Let $U^x_{n+1}=\bigcup\{U^y_0:y\in U^x_n\}$ and
$V_x=\bigcup\{U^x_n:n\in\omega\}.$ Notice that by this construction,
for every $x\in X$ and every $C\in\mathcal{C}$ if $C\cap
V_x\neq\varnothing$ then $C\subset V_x.$

Now if $x\in V_y$ then, there is $n\in\omega$ such that $x\in U^y_n$
which implies that $U^x_0\subset U^y_{n+1}.$ It follows that
$V_x\subset V_y.$ On the other hand, by the construction, $x\in
U^y_n$ also implies that $y\in U^x_n$ and thus $V_y\subset V_x.$ We
can conclude that if $V_y\cap V_x\neq\varnothing$ then $V_y=V_x.$

We have shown that the cover $\{V_x:x\in X\}$ induces a partition
$\mathcal{P}=\{P_\alpha:\alpha\in I\}$ of the space $X$ by countable
open subsets of $X.$ We can write
$P_\alpha=\{x^\alpha_n:n\in\omega\}$ and define
$D_n=\{x^\alpha_n:\alpha\in I\}$ for every $n\in\omega.$ It is clear
that $X=\bigcup\{D_n:n\in\omega\}$ and that each $D_n$ is discrete
because $P_\alpha$ isolates $x^\alpha_n$ in $D_n.$

\end{demos}

We have already observed that Eberlein-Grothendieck spaces of
cardinality $\omega_1$ are hereditarily meta-Lindelöf. This
observation allows us to obtain the following corollaries.

\begin{corol}

Every Eberlein-Grothendieck locally compact scattered space of
cardinality $\omega_1$ and height lower than $\omega_1\cdot\omega_1$
is $\sigma$-discrete.

\end{corol}

\begin{corol}
Every Eberlein-Grothendieck locally countable scattered space of
cardinality $\omega_1$ is $\sigma$-discrete.
\end{corol}

\begin{corol}
Every Eberlein-Grothendieck right-separated transfinite sequence
$X=\{x_\alpha:\alpha<\lambda\}$ with $\lambda<\omega_1\cdot\omega_1$
is $\sigma$-discrete.
\end{corol}

If we add separability as a hypothesis in Corollary 3.7, then we
have the following direct consequence of [Sp, Lemma 2.14].

\begin{corol}
Suppose that $K$ is a compact space of weight $\omega_1,$ if
$X\subset C_p(K)$ is locally countable separable and scattered then
$X$ is countable.
\end{corol}

As we already mentioned in the Introduction, when $K$ is a compact
scattered space, the $\sigma$-fragmentability of $C_p(K)$ is
equivalent to $C_p(K,\{0,1\})$ being $\sigma$-scattered, and
analogously the SLD property of $C_p(K)$ is equivalent to
$C_p(K,\{0,1\})$ being $\sigma$-discrete, cf. [Mtz].  The following
definition is a slight variation of [Fa, Definition 5.1.1]

\begin{defis}

Let $X$ be a space with a metric $\rho.$ An increasing family
$\{U_\beta: \beta<\alpha\}\subset \tau^*(X)$ that covers $X$ and
$diam_\rho(U_\gamma\setminus\displaystyle\bigcup_{\beta<\gamma}
U_\beta)<\varepsilon$ is called an $\varepsilon$-open partitioning
of length $\alpha$ of the space $X.$

\end{defis}

\begin{corol}
If $K$ is a scattered compact space of weight $\omega_1$ then
$C_p(K)$ is SLD if and only if for every $\varepsilon>0$
there exists a family $\{X_n:n\in\omega\}$ that covers
$C_p(K,\{0,1\})$ and for each $n\in\omega$ the set $X_n$ has an
$\varepsilon$-open partitioning of length $\omega_1.$
\end{corol}

\begin{demos}
Take $\varepsilon<1,$ there exists a family $\{X_n:n\in\omega\}$
that covers $C_p(K)$ and for each $n\in\omega$ the set $X_n$ has a
$\varepsilon$-open partitioning $\{U^n_\alpha: \alpha<\omega_1\}.$
For every $\alpha<\omega_1$ we have $|U^n_\alpha\cap
C_p(K,\{0,1\})|\leq1$ which implies that for each $n\in\omega,$ the
space $X_n\cap C_p(K,\{0,1\})$ is a scattered locally countable
space thus $\sigma$-discrete by Corollary 3.8 and hence
$C_p(K,\{0,1\})$ is $\sigma$-discrete as well. Apply [Mtz, Theorem
6] to conclude that $C_p(K)$ is SLD.

Assume that $C_p(K)$ is SLD and take $\varepsilon>0$ and a
countable cover $\{X_n:n\in\omega\}$ such that for each $n\in\omega$
the space $X_n$ has an $\varepsilon$-small open cover. For every
$n\in\omega$ we have $l(X_n)\leq w(C_u(K))=w(K)=\omega_1,$ therefore
$X_n$ has an $\varepsilon$-small open cover $\mathcal{V}_n$ of
cardinality $\omega_1.$  Let $\mathcal{V}_n=\{V^n_\beta:
\beta<\omega_1\}$ and define
$U^n_\beta=\displaystyle\bigcup_{\alpha\leq\beta}V^n_\alpha.$ The
family $\{U^n_\beta: \beta<\omega_1\}$ is an $\varepsilon$-open
partitioning of length $\omega_1.$

\end{demos}

\section{Topological games on Eberlein-Grothendieck spaces}

In [Ha] R.W. Hansell proved, among other things, that for every
compact $K$ that is a continuous image of a Valdivia compact, the
space $C_p(K)$ is hereditarily weakly $\theta$-refinable. Thus if
$K$ is a continuous image of a Valdivia compact then each scattered
subspace of $C_p(K)$ is $\sigma$-discrete. This implies in
particular that if $X$ is an Eberlein-Grothendieck Lindelöf $\Sigma$
scattered space then it is $\sigma$-dicrete. Indeed, suppose that
$X\subset C_p(Q)$ for some compact $Q$ and let $\varphi$ be the
diagonal product of the elements of $X$. If $K=\varphi(Q)$ then the
dual map $\varphi^*$ embeds the space $C_p(K)$ into $C_p(Q)$ as a
closed subspace that contains $X$ (see [Tk2, Problem 163]). It
follows that $C_p(K)$ contains a homeomorphic copy of $X$ that is a
Lindelöf $\Sigma$ subspace of $C_p(K)$ which separates the points of
$K.$ This implies that $C_p(K)$ is Lindelöf $\Sigma$ by [Ar,
Corollary IV.2.10] and therefore $K$ is a Gul'ko (and hence
Valdivia) compact space; thus $X$ is hereditarily weakly
$\theta$-refinable and consequently $\sigma$-discrete. As a
consequence, every Eberlein-Grothendieck K-analytic or cosmic
scattered space is $\sigma$-discrete.

Since every Lindelöf \v{C}ech-complete space $X$ is Lindelöf
$\Sigma$ (see for example [Tk1, Theorem 1]) we already know that if
$X$ is scattered then it is $\sigma$-discrete. However, we can say
more in this case; in this section we will apply the methods
developed in [Te] to prove that Lindelöf \v{C}ech-complete scattered
spaces are in fact $\sigma$-compact.

\begin{defis}
On a Tychonoff space $Y$, consider a family $\mathcal{C}$ of subsets
of $Y.$ We define the game $\mathcal{G}(\mathcal{C},Y)$ of two
players $I$ and $II$ who take turns in the following way: at the
move number $n,$ Player $I$ chooses $C_n\in\mathcal{C}$ and Player
$II$ chooses a set $U_n\in\tau(C_n,Y).$ The game ends after the
$n$-th move of each player has been made for every $n\in\omega$ and
Player $I$ wins if $X=\bigcup\{U_n: n\in\omega\};$ otherwise the
winner is Player $II$.
\end{defis}

\begin{defis}
A strategy $t$ for the first player in the game
$\mathcal{G}(\mathcal{C},Y)$ on a space $X$ is defined inductively
in the following way. First the set
$t(\varnothing)=F_0\in\mathcal{C}$ is chosen. An open set
$U_0\in\tau(X)$ is legal if $F_0\subset U_0.$ For every legal set
$U_0$ the set $t(U_0)=F_1\in\mathcal{C}$ has to be defined. Let us
assume that legal sequences $(U_0,...,U_i)$ and sets
$t(U_0,...,U_i)$ have been defined for each $i\leq n.$ The sequence
$(U_0,...,U_{n+1})$ is legal if so is the sequence $(U_0,...,U_i)$
for each $i\leq n$ and $F_{n+1}=t(U_0,...,U_n)\subset U_{n+1}.$ A
strategy $t$ for Player $I$ is winning if it ensures victory for $I$
in every play it is used.
\end{defis}

Recall that given a space $X$ and $Y\subset X,$ the set $Y$ has
countable character in $X$ if there exists a countable family
$\mathcal{U}\subset\tau(Y,X)$ such that for every $V\in\tau(Y,X)$
there exists $U\in\mathcal{U}$ such that $V\subset U.$ A space $X$
is of countable type if every compact subset of $X$ is contained in
a compact space of countable character in $X.$ Notice that in a
space of countable type $X$ not necessarily every compact subset of
$X$ has countable character in $X,$ therefore the following
proposition does not follow form the results in [Te], however we
will apply the ideas in [Te] to prove it.

\begin{prop}
If a space $X$ has countable type and $\mathcal{K}$ is the family of
its compact subspaces, then the first player has a winning strategy
for the game $\mathcal{G}(\mathcal{K},X)$ if and only if $X$ is
$\sigma$-compact.
\end{prop}

\begin{demos}
It suffices to prove necessity. Suppose that the space $X$ has
countable type and that the first player has a winning strategy for
the game $\mathcal{G}(\mathcal{K},X).$ Let $s$ be a winning strategy
for the first player in the game $\mathcal{G}(\mathcal{K},X).$ For
every $F\in\mathcal{K}$ there is $K(F)\in\mathcal{K}$ such that
$F\in K(F)$ and $\chi(K(F))\leq\omega$ hence, for each
$F\in\mathcal{K}$ it is possible to find a countable family
$\mathcal{U}_F=\left\{U_n^F: n<\omega \right\}\subset \tau (K(F),X)$
such that $\bigcap \mathcal{U}_F=K(F)$. For the compact set
$F_0=s(\varnothing)$ define $A_0=\{K(F_0)\}$ and for every
$n\in\mathbb{N}$ define
$A_n=\{K(s(U_{l_0}^{{F_0}},\ldots,U_{l_{n-1}}^{F_{n-1}})):
l_0,\ldots,l_{n-1}<\omega, F_i\in A_i$ for each $ i<n $ and
$((F_0,U_{l_0}^{F_0}),...,(F_{n-1},U_{l_{n-1}}^{F_{n-1}}))$ is an
initial segment of a match of the game $\mathcal{G}(\mathcal{K},X)$
in which the first player applies the strategy $s\}.$ Observe that
$\left|A_n\right|\leq\omega$ for every $n\in \omega.$ Indeed,
$\left|A_0\right|\leq\omega.$ If we assume that
$\left|A_n\right|\leq\omega$ then we have that
$\left|A_{n+1}\right|\leq \left|A_n\right|\cdot
\left|\mathcal{U}_{F_{n}}\right|\leq\omega^2\leq\omega$. Therefore
$\left|\bigcup\{A_n:n\in\omega\}\right|\leq\omega$.

For every $n\in\omega$ define $B_n=\bigcup A_n$ and $B=\bigcup B_n.$
We will show that $X=B$. Suppose that $y\in X\diagdown B;$ this
implies $y\neq K(F_0)$ thus there is $U_0\in\mathcal{ U}_{F_0}$ such
that $y\notin U_0.$ Let $F_1=s(U_0).$ The set $K(F_1)\in A_1$
therefore $y\notin K(F_1)$ and there is $U_1\in\mathcal{U}_{F_1}$
such that $y\notin U_1.$ Let $F_2=s(U_1,U_2)$. Suppose $F_k,
U_{k-1}$ have been defined by this procedure in such a way that
$y\notin U_j$ with $j=1,\ldots,k-1$, it is then possible to find
$U_k\in\mathcal{U}_{F_k}$ such that $y\notin U_k$. By the definition
of $\{U_n:n\in\omega\}$ we have that $\mathcal{P}=\{(F_n,U_n):n\in
\omega\}$ is a match of $\mathcal{G}(\mathcal{K},X)$ in which the
first player applies $s$, but $y\notin U_n$ for every $n\in \omega.$
This contradiction shows that $X$ is $\sigma$-compact.

\end{demos}

Let $\mathcal{S}$ be the family of singletons of a space Lindelöf
scattered space $X.$ In [Te, Theorem 9.3] Telgarsky proved that the
first player has a winning strategy for the game
$\mathcal{G}(\mathcal{S},X).$ Thus we have the following corollary.

\begin{corol}

Every Lindelöf scattered space of countable type is
$\sigma$-compact.

\end{corol}

The most important class of countable type spaces are the
\v{C}ech-complete spaces, so we obtain the following theorem.

\begin{theos}

Every Lindelöf \v{C}ech-complete scattered space is
$\sigma$-compact.

\end{theos}

Let $\mathcal{S}$ be the family of singletons of a compact space
$K.$ In [Te] it is proved, among other things, that the first player
has a winning strategy for the game $\mathcal{G}(\mathcal{S},K)$ if
and only if $K$ is scattered. We can apply this fact in the context
of Eberlein-Grothendieck spaces to obtain the following corollary.

\begin{corol}

Let $\mathcal{S}$ be the family of all singletons of an
Eberlein-Grothendieck \v{C}ech-complete space $X.$ If the first
player has a winning strategy for the game
$\mathcal{G}(\mathcal{S},X)$ then $X$ is $\sigma$-discrete.

\end{corol}

\begin{demos}
The space $X$ being \v{C}ech-complete has countable type. Besides by
[Te, Corollary 2.2], the first player also has a winning strategy
for the game $\mathcal{G}(\mathcal{K},X),$ where $\mathcal{K}$ is
the family of compact subsets of $X.$ Apply Theorem 4.3 to conclude
that $X$ is $\sigma$-compact. We can write
$X=\displaystyle\bigcup_{n\in\omega}K_n$ where each $K_n$ is
Eberlein compact. It is easy to see that for each $n\in\omega$ the
first player has a winning strategy for the game
$\mathcal{G}(\mathcal{S}_n,K_n),$ where $\mathcal{S}_n$ is the
family of singletons of $K_n$ implying that each $K_n$ is scattered
and consequently $\sigma$-discrete and so is $X.$

\end{demos}

\section{Example}

In Section 3 we applied the hereditary metalindelöfness of certain
scattered spaces to prove they are $\sigma$-discrete. It is not very
clear that this property implies $\sigma$-discreteness of scattered
Tychonoff spaces. This is not true for general spaces as we can
deduce from the following example.

\begin{example}
There exists a scattered space of class $T_1$ which is hereditarily
meta-Lindelöf that is not $\sigma$-discrete.

\end{example}

\begin{demos}
We will define a $T_1$ hereditarily meta-Lindelöf topology on the
ordinal $\omega_1\cdot\omega_1$. First we will define some auxiliary
sets. Fix an ordinal $\gamma<\omega_1$ and define
$V_{(\alpha,\gamma)}$ as follows:

\begin{itemize}

\item $V_{(\alpha,\gamma)}=\varnothing$ if $\alpha<\gamma.$

\item $V_{(\gamma,\gamma)}=[\omega_1\cdot\gamma,\omega_1\cdot(\gamma+1)).$

\item $V_{(\gamma+\beta,\gamma)}=[\omega_1\cdot\gamma+\beta,\omega_1\cdot(\gamma+1))$ for $\beta<\omega_1.$

\end{itemize}

Define
$V_\alpha=\displaystyle\bigcup_{\gamma<\omega_1}V_{(\alpha,\gamma)}.$

Observe that the family $\{V_\alpha:\alpha<\omega_1\}$ thus defined
is point countable. Indeed, let $\xi:\in\omega_1\cdot\omega_1.$
There are countable ordinals $\varsigma,\varrho$ such that
$\xi=\omega_1\cdot\varsigma+\varrho.$ It is clear that $\xi\notin
V_\alpha$ for $\alpha>(\varsigma+\varrho)+1.$

For every $\lambda<\omega_1\cdot\omega_1$ let
$U_\lambda=[0,\lambda)$ and $W^\lambda_\alpha=V_\alpha\cap
U_\gamma.$ It is easy to see that the family
$\mathcal{W}=\{W^\lambda_\alpha: \alpha<\omega_1,
\lambda<\omega_1\cdot\omega_1\}$ is a base for a topology $\tau$ on
$\omega_1\cdot\omega_1.$ Note that by this construction
$U_\lambda\in\tau$ for each $\lambda\in\omega_1\cdot\omega_1;$ hence
the space $(\omega_1\cdot\omega_1, \tau)$ with the order of
$\omega_1\cdot\omega_1$ is right-separated, and therefore this space
is scattered.

To verify that the space $(\omega_1\cdot\omega_1, \tau)$ is
hereditarily meta-Lindelöf it suffices to show that for any family
$\mathcal{C}\subset\mathcal{W}$ there is a point countable family
$\mathcal{O}\subset\tau$ such that
$\bigcup\mathcal{C}=\bigcup\mathcal{O}.$

Take a family of basic open sets $\mathcal{C}=\{W^\lambda_\alpha:
(\alpha,\lambda)\in I\}\subset\mathcal{W}.$ For each
$\alpha<\omega_1$ define $C_\alpha=\bigcup\mathcal{C}\cap V_\alpha$
and $J_\alpha=\{\lambda\in\omega_1\cdot\omega_1:
W^\lambda_\alpha\in\mathcal{C}\}.$ Notice that the family
$\{C_\alpha:\alpha<\omega_1\}$ is point countable, therefore it will
be enough to show a point countable family $\mathcal{O}_\alpha$ such
that $C_\alpha=\bigcup\mathcal{O}_\alpha$ for each
$\alpha<\omega_1.$

There are three possible mutually exclusive cases:

\begin{enumerate}
\item If there is $\xi\in J_\alpha$ such that $W^\lambda_\alpha\subset
W^\xi_\alpha$ for every $\lambda\in J_\alpha$ then define
$\mathcal{O}_\alpha=\{ W^\xi_\alpha\}$

\item If $J_\alpha$ is bounded in $V_\alpha$ and
$cf(J_\alpha)=\omega,$ then let $J'_\alpha$ be a countable cofinal
subset of $J_\alpha$ and $\mathcal{O}_\alpha=\{W^\lambda_\alpha:
\lambda\in J'_\alpha\}.$

\item If the cofinality of the set $J_\alpha$ is $\omega_1$ then it is
possible to find an increasing $\omega_1$-sequence of ordinals
$\{\lambda_\eta:
\eta<\omega_1\}\subset(\omega_1\cdot\beta,\omega_1\cdot(\beta+1)),$
with $\beta+1\leq\alpha+1,$ that is cofinal in $J_\alpha.$ For each
$\eta<\omega_1$ there is $\delta(\eta)<\omega_1$ such that
$\lambda_\eta=\omega_1\cdot\beta+\delta(\eta).$ Let
$\mathcal{O}_\eta=\{W^{\lambda_{(\eta+1)}}_\varsigma:
\beta+\delta(\eta)\leq\varsigma<\beta+\delta(\eta+1)\}.$ Define
$\mathcal{O}_\alpha=\displaystyle\bigcup_{\eta<\omega_1}\mathcal{O}_\eta.$
To see that in this case the family $\mathcal{O}_\alpha$ is point
countable, take $\rho\in\bigcup\mathcal{O}_\alpha.$ Since each
$\mathcal{O}_\eta$ is countable, it suffices to show that the set
$\{\eta<\omega_1:\rho\in\bigcup\mathcal{O}_\eta\}$ is countable. We
can find countable ordinals $\zeta$ and $\xi$ such that
$\rho=\omega_1\cdot\zeta+\xi.$ Observe that by the definition of
each $\delta(\eta)$ the family $\{\delta(\eta):\eta<\omega_1\}$ is
increasing and cofinal in $\omega_1.$ Thus there is $\eta<\omega_1$
such that $\zeta+\xi+1<\delta(\eta).$ It is clear that
$\rho\notin\mathcal{O}_{\eta'}$ for every $\eta'>\eta.$

\end{enumerate}

To see that $X$ is not $\sigma$-discrete, suppose that
$X=\displaystyle{\bigcup_{n\in\omega}D_n}.$ For each
$\alpha<\omega_1$ define $D^\alpha_n=D_n\cap[\omega_1\cdot\alpha,
\omega_1\cdot(\alpha+1)).$ For every $\alpha\in\omega_1$ we can find
$\Phi(\alpha)\in\omega$ such that
$|D^\alpha_{\Phi(\alpha)}|=\omega_1.$ Thus, there is $m\in\omega$
such that
$|\{\alpha<\omega_1:|D^\alpha_m|=\omega_1\}|=|\Phi^{-1}(m)|=\omega_1.$
It is possible to find ordinals $\beta\in\omega_1$ y $\gamma\in D_m$
such that $|\{\alpha\in\beta: |D^\alpha_m|=\omega_1\}|=\omega$ and
$\omega_1\cdot\beta<\gamma\in D_m.$

If $\gamma\in W^\xi_\delta\in\tau$ then for every $\alpha<\beta $
such that $|D^\alpha_m|=\omega_1$ there exist $\lambda<\omega_1$
such that $W^\xi_\delta\cap[\omega_1\cdot\alpha,
\omega_1\cdot(\alpha+1))=[\omega_1\cdot\alpha+\lambda,
\omega_1\cdot(\alpha+1)).$ Since $|D^\alpha_m|=\omega_1,$ it follows
that $D_m\cap[\omega_1\cdot\alpha+\lambda,
\omega_1\cdot(\alpha+1))\neq\varnothing.$ This shows that the set
$D_m$ is not a discrete subspace of $(\omega_1\cdot\omega_1, \tau).$

\end{demos}

Example 6.1 shows that hereditary metalindelöfness alone is not
enough to guarantee $\sigma$-discreteness of a scattered topological
spaces in general. Unfortunately, this example does not apply to the
context of Eberlein-Grothendieck spaces because it is not even a
Hausdorff space. Indeed, any two uncountable ordinals cannot be
separated by disjoint open subsets. Hence a natural question arises:
is every Eberlein-Grothendieck hereditarily metalindelöf scattered
space $\sigma$-discrete? (see Problem 6.4).

\section{Open Problems}

It is already known that Eberlein scattered spaces are
$\sigma$-discrete, the next evident step would be to consider the
case of Eberlein-Grothendieck Lindelöf spaces. A positive answer
would follow if $C_p(K)$ was hereditarily weakly $\theta$-refinable
whenever $K$ is a compact subspace of $C_p(L)$ for a scattered
Lindelöf space $L.$

\begin{problem}
Is every Eberlein-Grothendick Lindelöf scattered space
$\sigma$-discrete?
\end{problem}

\begin{problem}
Suppose that $L$ is a Lindelöf scattered space and $K$ is a compact
subspace of $C_p(L).$ Is  $C_p(K)$ hereditarily weakly $\theta$
refinable?
\end{problem}

In Section 4 we showed that Problem 1.1 has a positive partial
answer for the case of Lindelöf \v{C}ech-complete scattered spaces.
What if we remove the hypothesis that the space is Lindelöf?

\begin{problem}
Is every Eberlein-Grothendieck \v{C}ech-complete scattered space
$\sigma$-discrete?
\end{problem}

As noticed in Section 5, if we consider spaces that are not
Tychonoff, hereditarily metalindelöfness does not neccesarily imply
$\sigma$-discreteness of scattered spaces in general. Does it in
Eberlein-Grothendieck spaces?

\begin{problem}
Is every Eberlein-Grothendieck hereditarily meta-Lindelöf scattered
space $\sigma$-discrete?
\end{problem}

Problem 6.3 would be a positive result, assuming MA, if the
following question had an affirmative answer.

\begin{problem}
Let $X$ be an Eberlein-Grothendieck hereditarily meta-Lindelöf
scattered space of height and cardinality equal to $\omega_1.$ For
every point $x\in X$ there is an open set $U_x$ that isolates $x$ in
its scattering level. There is a point countable open refinement
$\mathcal{V}$ of the cover $\{U_x: x \in X\}.$ Let $V_x$ be the
union of all the elements of $\mathcal{V}$ that contain $x.$ Define
the partially ordered set $P=\{ p \subset X: p $ is finite and
$V_x\cap p=\{x\}$ for every $x \in p \}$ and $q<p$ if $p\subset q.$
Has $P$ the countable chain condition?
\end{problem}

Take an Eberlein-Grothendieck right-separated transfinite sequence
$X=\{x_\alpha:\alpha<\lambda\}.$ In Section 3 we showed that $X$ is
hereditarily metalindelöf for $\lambda<\omega_2.$ Moreover we proved
that hereditarily metalindelöfness implies $\sigma$-discreteness of
$X$ for $\lambda<\omega_1\cdot\omega_1.$ It is not yet clear if
hereditarily metalindelöfness implies $\sigma$-discreteness of $X$
for $\lambda=\omega_1\cdot\omega_1.$

\begin{problem}
Suppose that $\omega_1\cdot\omega_1\leq\lambda<\omega_2$ and
$X=\{x_\alpha:\alpha<\lambda\}$ is an Eberlein-Grothendieck
right-separated transfinite sequence. Is $X$ $\sigma$-discrete?
\end{problem}


\begin{thebibliography}{Gue}

\bibitem[Al]{} K. Alster, \textit{Some remarks on Eberlein compacts.}
Fund. Math. \textbf{1:104}(1979), 43-46.

\bibitem[Ar]{}A.V. Arhangel'skii, \textit{Topological function
spaces,} Mathematics and its Applications (Soviet Series),
\textbf{78}, Kluwer Acad. Publ., Dordrecht, 1992.

\bibitem [Bu]{}
D. Burke, \textit{Covering properties,} Kunen y Vaughan (eds.),
Elsevier Science Publishers, Netherlands, 1984, 347-423.

\bibitem [En]{}R. Engelking, \textit{General Topology,} PWN,
Warszawa 1977.

\bibitem[DJP]{} A. Dow, H. Junnila, J. Pelant,
\textit{Weak covering properties of weak topologies,}
Proc. London Math. Soc., \textbf{3:75}(1997), 349-368.

\bibitem[Fa]{} M. J. Fabian,
\textit{Gâteaux differentiability of convex functions and topology.
Weak Asplund spaces,} Canadian Mathematical Society Series of
Monographs and Advanced Texts, A Wiley-Interscience Publication.
John Wiley \& Sons, Inc., New York, 1997.

\bibitem[Ha]{} R.W. Hansell, \textit{Descriptive sets and the
topology of nonseparable Banach spaces,} Serdica Math. J.,
\textbf{1:27}(2001), 1-66.

\bibitem[Hay]{} R. Haydon, \textit{Some problems about scattered
spaces,} Séminaire d'Initiation à l'Analyse, Exp. No. 9, 10 pp.,
Publ. Math. Univ. Pierre et Marie Curie, 95, Univ. Paris VI, Paris,
199?.

\bibitem[Ha]{} H. Z. Hdeib, C. M. Pareek,  \textit{A generalization
of scattered spaces,} Topology Proc. \textbf{1:14}(1989), 59-74.

\bibitem[JNR]{} J. E. Jayne, I. Namioka, C. A. Rogers,
\textit{$\sigma$-fragmentable Banach spaces,} Mathematika
\textbf{2:39}(1992) 197-215.

\bibitem[Mtz]{} J.F. Martínez, \textit{Sigma-fragmentability and
the property SLD in C(K) spaces,}  Topology Appl. \textbf{8:156}
(2009), 1505-1509.

\bibitem[MOTV]{} A. Moltó, J. Orihuela, S. Troyanski, M. Valdivia,
\textit{A nonlinear transfer technique for renorming}, Lecture Notes
in Mathematics 1951, Springer (2009).

\bibitem[Sp]{} S. Spadaro,  \textit{A note on discrete sets,}
Comment. Math. Univ. Carolin. \textbf{3:50}(2009), 463-475.

\bibitem[Te]{} R. Telgársky,
\textit{Spaces defined by topological games,} Fund. Math.,
\textbf{88:3}(1975), 193-223.

\bibitem [Tk1]{} V.V. Tkachuk, \textit{Lindelöf $\Sigma$-spaces: an
omnipresent class,} Rev. R. Acad. Cienc. Exactas Fís. Nat. Ser. A
Math. RACSAM \textbf{2:104}(2010), 221-244.


\bibitem [Tk2]{} V.V. Tkachuk, \textit{A $C_p$-Theory Problem Book,
Topological and Function Spaces} Springer, New York, 2011.


\bibitem[Ya]{}N.N. Yakovlev, \textit{On bicompacta in $\Sigma$-products
and related spaces,} Comment. Math. Univ. Carolinae,
\textbf{21:2}(1980), 263-283.

\end{thebibliography}
\end{document}